\documentclass[12pt]{amsart}
\usepackage{amssymb}

\usepackage{graphicx}
\usepackage{amscd}

\newtheorem{theorem}{Theorem}[section]
\newtheorem{corollary}[theorem]{Corollary}
\newtheorem{proposition}[theorem]{Proposition}

\input{tcilatex}

\begin{document}
\title[ Inequality involving packing dimension]{Multifractal formalism and inequality involving packing dimension}
\author{L. Ben Youssef.}
\address{Current address : L. Ben Youssef. ISCAE. University of Manouba. Tunisia.}
\email{Leila.BenYoussef@iscae.rnu.tn}
\date{June 7, 2008}

\begin{abstract}
This article fits in many studies of multifractal analysis of measure \cite
{F.B1, F.B2, F.B3, G.B, O.L1, O.L2, O.L3, J.P}. We took as a starting point
the work of F. Ben Nasr in \cite{F.B2} to give a new inequality involving $%
Dim(\overline{X}^{\alpha })$ which would be, in certain cases, finer than
the inequality
\begin{equation*}
Dim(\overline{X}^{\alpha })\leq \underset{q\geq 0}{\inf }(\alpha q+B_{\mu
}(q)),
\end{equation*}
established by L. Olsen in \cite{O.L1}. Besides we elaborated an application
of our result which gives a better inequality involving $Dim(\overline{X}%
^{\alpha })$. \newline
\newline
We are thankful to Mr F. Ben Nasr for the long and lucrative discussions
which we had during the development of this work.
\end{abstract}

\maketitle

\section{Multifractal formalism}

Let $\mu $\ be a Borel probability measure on $\mathbb{R}^{d}$. For $%
E\subset \mathbb{R}^{d}$, $q,$ $t\in \mathbb{R}$ and $\varepsilon >0$, by
adopting the convention
\begin{equation*}
\left\{
\begin{array}{l}
0^{q}=+\infty ,\text{ }q<0, \\
0^{0}=1,
\end{array}
\right.
\end{equation*}
put
\begin{equation*}
\overline{P}_{\mu ,\varepsilon }^{q,t}(E)=\sup \left\{ \underset{i}{\sum }%
\mu \left( B\left( x_{i},r_{i}\right) \right) ^{q}(2r_{i})^{t}\right\}
\end{equation*}
where the supremum is taken over all the centered $\varepsilon -$packing $%
(B\left( x_{i},r_{i}\right) )_{i\in I}$ of $E$.\newline
Also put
\begin{equation*}
\overline{P}_{\mu }^{q,t}(E)=\underset{\varepsilon \rightarrow 0}{\lim }%
\overline{P}_{\mu ,\varepsilon }^{q,t}(E).
\end{equation*}
Since $\overline{P}_{\mu }^{q,t}$ is a prepacking-measure, then we consider,

\begin{equation*}
P_{\mu }^{q,t}(E)=\underset{E\subset \left( \underset{i}{\cup }E_{i}\right)
}{\inf }\underset{i}{\sum }\overline{P}_{\mu }^{q,t}(E_{i}).
\end{equation*}
It is clear that
\begin{equation}
P_{\mu }^{q,t}(E)=\underset{E=\left( \underset{i}{\cup }E_{i}\right) }{\inf }%
\underset{i}{\sum }\overline{P}_{\mu }^{q,t}(E_{i})  \label{eq1}
\end{equation}
and
\begin{equation}
P_{\mu }^{q,t}(E)=\inf \left\{ \underset{i}{\sum }\overline{P}_{\mu
}^{q,t}(E_{i}):\left( \underset{i}{\cup }E_{i}\right) \ \text{is a partition
of }E\right\}  \label{eq2}
\end{equation}

The prepacking-measure $\overline{P}_{\mu }^{q,t}$ and the measure $P_{\mu
}^{q,t}$ assign respectively a dimension to each subset $E$. These
dimensions are respectively denoted by $\Delta _{\mu }^{q}(E)$ and $Dim_{\mu
}^{q}(E)$. They are respectively characterized by
\begin{equation*}
\overline{P}_{\mu }^{q,t}(E)=\left\{
\begin{array}{c}
\infty \text{ if }t<\Delta _{\mu }^{q}(E) \\
0\text{ if }t>\Delta _{\mu }^{q}(E)
\end{array}
\right.
\end{equation*}
and
\begin{equation*}
P_{\mu }^{q,t}(E)=\left\{
\begin{array}{c}
\infty \text{ if }t<Dim_{\mu }^{q}(E) \\
0\text{ if }t>Dim_{\mu }^{q}(E)
\end{array}
\right.
\end{equation*}
Note that L. Olsen established in \cite{O.L1} the following results
\begin{equation}
Dim_{\mu }^{q}(E)\leq \Delta _{\mu }^{q}(E)  \label{eq3}
\end{equation}
and
\begin{equation}
Dim_{\mu }^{q}(E)=\underset{E=\underset{n}{\cup }E_{n}}{\inf }\left\{
\underset{n}{\sup }Dim(E_{n})\right\} .  \label{eq4}
\end{equation}

The numbers $\Delta _{\mu }^{q}(E)$ and $Dim_{\mu }^{q}(E)$ are respectively
the multifractal extensions of the prepacking dimension $\Delta (E)$ and the
packing dimension $Dim(E)$ of $E$ (cf \cite{C.T}),\ in fact
\begin{equation*}
\Delta _{\mu }^{0}(E)=\Delta (E)\text{ \ \ \ \ \ \ and \ \ \ \ \ \ }Dim_{\mu
}^{0}(E)=Dim(E).
\end{equation*}

Write $\Lambda _{\mu }(q)=\Delta _{\mu }^{q}($supp$\mu )$ and $B_{\mu
}(q)=Dim_{\mu }^{q}($supp$\mu ).$ L. Olsen also established in \cite{O.L1}
the following results.

\begin{proposition}
\label{prop1} \ \ \ \ \newline
i. $B_{\mu }\leq \Lambda _{\mu }$, $B_{\mu }(1)=\Lambda _{\mu }(1)=0.$%
\newline
ii. $\Lambda _{\mu }(0)=\Delta ($supp$\mu )$ and $B_{\mu }(0)=Dim($supp$\mu
).$\newline
iii. The functions $\Lambda _{\mu }:q\mapsto \Lambda _{\mu }(q)$ and $B_{\mu
}:q\mapsto B_{\mu }(q)$ are convex and decreasing.
\end{proposition}

\begin{theorem}
\label{th1}\ For $\alpha \geq 0,$ put
\begin{equation*}
\overline{X}^{\alpha }=\left\{ x\in \text{supp}\mu :\underset{r\rightarrow 0%
}{\lim \sup }\frac{Log\mu (B(x,r))}{Logr}\leq \alpha \right\} .
\end{equation*}
If $\alpha q+B_{\mu }(q)\geq 0,$ then
\begin{equation*}
Dim(\overline{X}^{\alpha })\leq \underset{q\geq 0}{\inf }(\alpha q+B_{\mu
}(q)).
\end{equation*}
\end{theorem}

\section{An other inequality involving $Dim(\overline{X}^{\protect\alpha })$}

For all $\varepsilon >0$, let $(u_{\varepsilon })_{\varepsilon >0}$ be a
family of numbers such that $\varepsilon \leq u_{\varepsilon }$ and $%
\underset{\varepsilon \rightarrow 0}{\lim }u_{\varepsilon }=0$. Let $k\geq 1$
be an integer$.$ If $M$ $\subset $ supp$\mu $, for each centered $%
\varepsilon -$packing $(B(x_{i},r_{i}))_{i\in I}$ of $M$, we consider all
the families $\left( B(y_{i},\delta _{i})\right) _{i\in I}$ indexed by $I$
verifying the following property $(\mathcal{P}_{k})$ : \newline
there exists a finite partition of $I$ such that $I=I_{1}\cup ...\cup I_{s}$
with \ $1\leq s\leq k$ and $\left( B(y_{i},\delta _{i})\right) _{i\in I_{j}}$
a centered $u_{\varepsilon }-$packing of $M$ for all $1\leq j\leq s$.\newline
Then we define the quantity
\begin{equation*}
L_{\varepsilon ,(B(x_{i},r_{i}))_{i\in I}}^{k}(M)=\inf \left( \underset{i\in
I}{\sup }\left( \frac{Log\mu (B(y_{i},\delta _{i}))}{Log2r_{i}}\right)
\right)
\end{equation*}
where the infimum is taken over all the families verifying $(\mathcal{P}%
_{k}).$\newline
Now write
\begin{equation*}
L_{\varepsilon }^{k}(M)=\sup \left\{ L_{\varepsilon ,(B(x_{i},r_{i}))_{i\in
I}}^{k}(M)\right\}
\end{equation*}
where the supremum is taken over all the centered $\varepsilon -$packing $%
(B\left( x_{i},r_{i}\right) )_{i\in I}$ of $M.$ \newline
Remark that
\begin{equation}
L_{\varepsilon ,(B(x_{i},r_{i}))_{i\in I}}^{k}(M)\leq \underset{i\in I}{\sup
}\left( \frac{Log\mu (B(x_{i},r_{i}))}{Log2r_{i}}\right) .\   \label{eq5}
\end{equation}
On the other hand, when $\varepsilon <\varepsilon ^{\prime }$, $%
L_{\varepsilon ^{\prime }}^{k}(M)>L_{\varepsilon }^{k}(M)$, then we define
\begin{equation*}
L^{k}(M)=\underset{\varepsilon \rightarrow 0}{\lim }L_{\varepsilon }^{k}(M).
\end{equation*}
As the sequence $\left( L^{k}(M)\right) _{k}$ is decreasing, write
\begin{equation*}
L(M)=\underset{k\rightarrow +\infty }{\lim }L^{k}(M).
\end{equation*}

Thereafter, for $\eta >\alpha $ and $p\in \mathbb{N}\backslash \left\{
0\right\} $, write
\begin{equation*}
X_{\alpha }(\eta ,p)=\left\{ x\in \overline{X}^{\alpha }:2r\geq \frac{1}{p}%
\text{ or }(2r)^{\eta }\leq \mu (B(x,r))\right\} .
\end{equation*}
It is clear that $X_{\alpha }(\eta ,p)\subset X_{\alpha }(\eta ,p+1)$.
Besides, it follows from the equality
\begin{equation*}
\overline{X}^{\alpha }=\left\{ x\in \text{supp}\mu :\underset{r\rightarrow 0%
}{\lim \sup }\frac{Log\mu (B(x,r))}{Log2r}\leq \alpha \right\} ,
\end{equation*}
that
\begin{equation*}
\overline{X}^{\alpha }=\underset{p}{\cup }X_{\alpha }(\eta ,p).
\end{equation*}

\begin{proposition}
\label{prop2}\ For all $M\subset X_{\alpha }(\eta ,p),$%
\begin{equation*}
L(M)\leq \eta .
\end{equation*}
\end{proposition}

\begin{proof}
For $\varepsilon <\dfrac{1}{p}$ and $(B(x_{i},r_{i}))_{i\in I}$ a centered $%
\varepsilon -$packing of $M$, \newline we have for all $i\in I,$
\begin{equation*}
\frac{Log\mu (B(x_{i},r_{i}))}{Log2r_{i}}\leq \eta ,
\end{equation*}
hence
\begin{equation*}
\underset{i\in I}{\sup }\frac{Log\mu
(B(x_{i},r_{i}))}{Log2r_{i}}\leq \eta ,
\end{equation*}
from the inequality (\ref{eq5}), we deduce that
\begin{equation*}
L_{\varepsilon ,(B(x_{i},r_{i}))_{i\in I}}^{k}(M)\leq \eta ,
\end{equation*}
while considering the supremum over $I$, it results that
\begin{equation*}
L_{\varepsilon }^{k}(M)\leq \eta ,
\end{equation*}
letting $\varepsilon \rightarrow 0,$ we obtain
\begin{equation*}
L^{k}(M)\leq \eta ,
\end{equation*}
then letting $k\rightarrow +\infty ,$ it follows that
\begin{equation*}
L(M)\leq \eta .
\end{equation*}
\end{proof}

\begin{theorem}
Assume\label{th2} $s:=\underset{q}{\inf }B_{\mu }(q)<0.$ \newline
Put,
\begin{equation*}
T_{\mu }(\alpha ,\eta ,p)=\underset{M\subset X_{\alpha }(\eta ,p)}{\sup }L(M)%
\text{,}
\end{equation*}
\begin{equation*}
T_{\mu }(\alpha ,\eta )=\underset{p\rightarrow +\infty }{\lim }T_{\mu
}(\alpha ,\eta ,p),
\end{equation*}
\begin{equation*}
T_{\mu }(\alpha )=\underset{\eta \rightarrow \alpha ^{+}}{\lim }T_{\mu
}(\alpha ,\eta ),
\end{equation*}
then,
\begin{equation*}
Dim(\overline{X}^{\alpha })\leq \frac{1}{\alpha }T_{\mu }(\alpha )\text{ }%
\underset{q\geq 1}{\inf }\left( \alpha q+B_{\mu }(q)\right) .
\end{equation*}
\end{theorem}

The limits $T_{\mu }(\alpha ,\eta )$ and $T_{\mu }(\alpha )$ are well
defined, in fact the sequence $\left( T_{\mu }(\alpha ,\eta ,p)\right)
_{p\geq 1}$ is increasing, since $X_{\alpha }(\eta ,p)\subset X_{\alpha
}(\eta ,p+1)$ and for all $\eta <\eta ^{\prime }$, $X_{\alpha }(\eta
,p)\subset X_{\alpha }(\eta ^{\prime },p)$, thus the quantity is $T_{\mu
}(\alpha ,\eta )$ decreasing when $\eta \rightarrow \alpha .$

Let's note that, according to the Proposition \ref{prop2}, we still have
\newline
$\dfrac{1}{\alpha }T_{\mu }(\alpha )\leq 1$, \thinspace what would permit in
some cases, in comparison with the Theorem \ref{th1}\ established by L.
Olsen in \cite{O.L1}, to compensate the loss achieved on the quantity $%
\underset{q\geq 1}{\inf }\left( \alpha q+B_{\mu }(q)\right) .$ In fact the
new inequality involving $Dim(\overline{X}^{\alpha })$ is better in
particular when
\begin{equation*}
\frac{1}{\alpha }T_{\mu }(\alpha )<1\text{ and }\underset{q\geq 1}{\inf }%
\left( \alpha q+B_{\mu }(q)\right) =\underset{q\geq 0}{\inf }\left( \alpha
q+B_{\mu }(q)\right) .
\end{equation*}
We will develop in the following paragraph an example where these conditions
will be verified.

Before proving the Theorem \ref{th2}, we establish the following proposition.

\begin{proposition}
\label{prop3}\ For $z>s$, put $\psi (z)=\inf B_{\mu }^{-1}(\left] -\infty ,z%
\right[ )$. Then
\begin{equation*}
\inf \left\{ \psi (\eta t)+t:\frac{s}{\eta }<t<0\right\} =\frac{1}{\eta }%
\text{ }\underset{q\geq 1}{\inf }\left( \eta q+B_{\mu }(q)\right) .
\end{equation*}
\end{proposition}

\begin{proof}
As $B_{\mu }$ is convex, decreasing on $\left[ 0,+\infty \right[ $
and taking strictly negative values for $\dfrac{s}{\eta }<t<0$,
there exists an unique $q>1$ such that $\eta t=B_{\mu }(q).$\newline
We deduce that for all $n\in \mathbb{N}\backslash \left\{ 0\right\} ,$ $%
B_{\mu }(q+\dfrac{1}{n})<B_{\mu }(q)=\eta t$.\newline Thus
$q+\dfrac{1}{n}\in B_{\mu }^{-1}(\left] -\infty ,\eta t\right[ )$
therefore $\psi (\eta t)\leq q,$ \newline On the other hand, from
the equalities
\begin{equation*}
\psi (\eta t)=\inf B_{\mu }^{-1}(\left] -\infty ,\eta t\right[
)=\inf \left\{ \theta :B_{\mu }(\theta )<\eta t\right\} =\inf
\left\{ \theta :\theta >q\right\} ,
\end{equation*}
it follows that $\psi (\eta t)\geq q.$ So $\psi (\eta t)=q.$\newline
There are two possible cases :\newline If $s=-\infty $, then it is
clear that
\begin{equation*}
\inf \left\{ \psi (\eta t)+t:\frac{s}{\eta }<t<0\right\} =\inf \left\{ q+%
\frac{B_{\mu }(q)}{\eta }:q>1\right\} .
\end{equation*}
If $s>-\infty $, put $q_{s}=\underset{q>1}{\inf }\left\{ q:B_{\mu
}(q)=s\right\} $. As $B_{\mu }$ is convex, it follows that for all
$q\geq
q_{s}$, $B_{\mu }(q)=s$, then $q+\dfrac{B_{\mu }(q)}{\eta }\geq q_{s}+\dfrac{%
s}{\eta }.$ It results that
\begin{equation*}
\inf \left\{ q+\frac{B_{\mu }(q)}{\eta }:q>1\right\} =\inf \left\{ q+\frac{%
B_{\mu }(q)}{\eta }:1<q\leq q_{s}\right\} .
\end{equation*}
Consider a sequence $(q_{n})$ such that $q_{n}\rightarrow q_{s}$ and $%
1<q_{n}<q_{s},$ as $B_{\mu }$ is continuous, surely we obtain $q_{n}+\dfrac{%
B_{\mu }(q_{n})}{\eta }\rightarrow q_{s}+\dfrac{B_{\mu
}(q_{s})}{\eta }.$ It follows that
\begin{equation*}
\inf \left\{ q+\frac{B_{\mu }(q)}{\eta }:q>1\right\} =\inf \left\{ q+\frac{%
B_{\mu }(q)}{\eta }:1<q<q_{s}\right\} ,
\end{equation*}
i.e.
\begin{equation*}
\inf \left\{ q+\frac{B_{\mu }(q)}{\eta }:q>1\right\} =\inf \left\{
\psi (\eta t)+t:\frac{s}{\eta }<t<0\right\} .
\end{equation*}
Otherwise, $(1+\dfrac{1}{n})+\dfrac{B_{\mu }(1+\frac{1}{n})}{\eta }%
\rightarrow 1+\dfrac{B_{\mu }(1)}{\eta },$ so

\begin{equation*}
\inf \left\{ q+\frac{B_{\mu }(q)}{\eta }:q>1\right\} =\inf \left\{ q+\frac{%
B_{\mu }(q)}{\eta }:q\geq 1\right\} ,
\end{equation*}
finally,
\begin{equation*}
\inf \left\{ \psi (\eta t)+t:\frac{s}{\eta }<t<0\right\} =\frac{1}{\eta }%
\text{ }\underset{q\geq 1}{\inf }\left( \eta q+B_{\mu }(q)\right) .
\end{equation*}
\end{proof}

\begin{proof}[Proof of the Theorem 2.2]
We stand in the interesting case where $\overline{X}^{\alpha }\neq
\varnothing $, it follows that $\alpha q+B_{\mu }(q)\geq 0,$ for all
$q\geq 0 $. Thus for $\eta >\alpha $, $\eta q+B_{\mu }(q)\geq 0$.
Then, from the
Proposition \ref{prop3}, it follows that if $\frac{s}{\eta }<t<0,$ then $%
\psi (\eta t)+t\geq 0.$

For $\gamma >0$ and $\dfrac{s}{\eta }<t<0,$ if $\gamma >\psi (\eta
t)+t,$ then $B_{\mu }(\gamma -t)<\eta t.$ It results that $P_{\mu
}^{\gamma -t,\eta t}($supp$\mu )=0,$ then $P_{\mu }^{\gamma -t,\eta
t}(X_{\alpha }(\eta ,p))=0. $ According to the equality (\ref{eq1}),
we can write $X_{\alpha
}(\eta ,p)=\underset{u\in U}{\cup }M_{u}$ such that for all $u\in U$ and $%
\lambda >T_{\mu }(\alpha ,\eta ,p),\overline{P}_{\mu }^{\gamma
-t,\eta
t}(M_{u})<\infty .$ First of all let's prove that for all $u\in U$, $%
\triangle (M_{u})\leq \gamma \lambda .$ As $M_{u}\subset X_{\alpha
}(\eta
,p) $ and $\lambda >L(M_{u}),$ then there exist an integer $k\geq 1$ and $%
\varepsilon _{0}<\frac{1}{p}$ such that for all $\varepsilon
<\varepsilon _{0},$
\begin{equation*}
L_{\varepsilon }^{k}(M_{u})<\lambda \text{ and }\overline{P}_{\mu
,\varepsilon }^{\gamma -t,\eta t}(M_{u})<\infty .
\end{equation*}
For all $(B(x_{i},r_{i}))$ centered $\varepsilon -$packing of
$M_{u},$ there
exists a family $\left( B(y_{i},\delta _{i})\right) _{i\in I}$ such that $%
I=I_{1}\cup ...\cup I_{s}$ with \ $1\leq s\leq k$ and $\left(
B(y_{i},\delta _{i})\right) _{i\in I_{j}}$ a centered
$u_{\varepsilon }-$packing of $M_{u}$ for all $1\leq j\leq s$ and
for all $i\in I,$ $\dfrac{Log\mu (B(y_{i},\delta
_{i}))}{Log2r_{i}}<\lambda .$ It follows that $\mu (B(y_{i},\delta
_{i}))>(2r_{i})^{\lambda }$ and $\left( 2\delta _{i}\right) ^{\eta
}<\mu (B(y_{i},\delta _{i})).$ Thus for $\gamma >0$ and
$\dfrac{s}{\eta }<t<0,$ we obtain
\begin{equation*}
(2r_{i})^{\gamma \lambda }\leq \mu (B(y_{i},\delta _{i}))^{\gamma
-t}\left( 2\delta _{i}\right) ^{\eta t},
\end{equation*}
then,
\begin{equation*}
\underset{i\in I}{\sum }(2r_{i})^{\gamma \lambda }\leq \underset{i\in I}{%
\sum }\mu (B(y_{i},\delta _{i}))^{\gamma -t}\left( 2\delta
_{i}\right) ^{\eta t}=\underset{j=1}{\overset{s}{\sum
}}\underset{i\in I_{j}}{\sum }\mu (B(y_{i},\delta _{i}))^{\gamma
-t}\left( 2\delta _{i}\right) ^{\eta t}
\end{equation*}
it results that
\begin{equation*}
\underset{i\in I}{\sum }(2r_{i})^{\gamma \lambda }\leq
k\overline{P}_{\mu ,\varepsilon }^{\gamma -t,\eta t}(M_{u})<\infty .
\end{equation*}
Then for all $u\in U$,
\begin{equation*}
\triangle (M_{u})\leq \gamma \lambda .
\end{equation*}
Therefore, from the inequality (\ref{eq3}),
\begin{equation*}
Dim(M_{u})<\gamma \lambda ,\text{ }u\in U.
\end{equation*}
And from the equality (\ref{eq4}), we deduce that
\begin{equation*}
Dim(X_{\alpha }(\eta ,p))<\gamma \lambda .
\end{equation*}
Thus
\begin{equation*}
Dim(X_{\alpha }(\eta ,p))\leq \gamma T_{\mu }(\alpha ,\eta ,p).
\end{equation*}
As $\overline{X}^{\alpha }=\underset{p\geq 1}{\cup }X_{\alpha }(\eta
,p)$, from the equality (\ref{eq4}), letting $p\rightarrow +\infty
,$ we obtain
\begin{equation*}
Dim(\overline{X}^{\alpha })\leq \gamma T_{\mu }(\alpha ,\eta ).
\end{equation*}
So for $\dfrac{s}{\eta }<t<0,$
\begin{equation*}
Dim(\overline{X}^{\alpha })\leq \inf \left\{ \psi (\eta t)+t:\frac{s}{\eta }%
<t<0\right\} T_{\mu }(\alpha ,\eta ).
\end{equation*}
Then, according to the Proposition \ref{prop3}, it results that\ for all $%
\eta >\alpha ,$
\begin{equation*}
Dim(\overline{X}^{\alpha })\leq \frac{1}{\eta }\text{ }\underset{q\geq 1}{%
\inf }\left( \eta q+B_{\mu }(q)\right) T_{\mu }(\alpha ,\eta ).
\end{equation*}
Finally, letting $\eta \rightarrow \alpha ,$ it follows that
\begin{equation*}
Dim(\overline{X}^{\alpha })\leq \frac{1}{\alpha }T_{\mu }(\alpha )\text{ }%
\underset{q\geq 1}{\inf }\left( \alpha q+B_{\mu }(q)\right) .
\end{equation*}
\end{proof}

\section{Example}

In this paragraph, we intend to construct a measure $\mu $ verifying the
following conditions
\begin{equation*}
\frac{1}{\alpha }T_{\mu }(\alpha )<1\text{ and }\underset{q\geq 1}{\inf }%
\left( \alpha q+B_{\mu }(q)\right) =\underset{q\geq 0}{\inf }\left( \alpha
q+B_{\mu }(q)\right) .
\end{equation*}
What will permit, thanks to the Theorem \ref{th2}, to establish that
\begin{equation*}
Dim(\overline{X}^{\alpha })\leq \frac{1}{\alpha }T_{\mu }(\alpha )\text{ }%
\underset{q\geq 1}{\inf }\left( \alpha q+B_{\mu }(q)\right) <\underset{q\geq
0}{\inf }(\alpha q+B_{\mu }(q)).
\end{equation*}

Put $\mathcal{A}$ the set of the words constructed with $\left\{ 0,1\right\}
$ as alphabet. The length of a word $j$ is denoted by $\left| j\right| .$
For all $j\in \mathcal{A},$ put $N_{0}(j)$ the number of times the letter $0$
appears in $j.$ If $j,$ $j^{\prime }\in \mathcal{A}$, write $jj^{\prime }$
the word starting by $j$ and gotten while putting $j^{\prime }$ after $j$.
For all $j\in \mathcal{A}$ such that $j=j_{1}j_{2}...j_{n}$, put $I_{j}$ the
diadic interval of order $n$ defined by
\begin{equation*}
I_{j}=\left[ \underset{k=1}{\overset{n}{\sum }}\frac{j_{k}}{2^{n}},\underset{%
k=1}{\overset{n}{\sum }}\frac{j_{k}}{2^{n}}+\frac{1}{2^{n}}\right[ .
\end{equation*}
We denote by $\mathcal{F}_{n}$ the family of all the diadic intervals of
order $n$ and for all $x\in \left[ 0,1\right[ $ we call $I_{n}(x)$ the
element of $\mathcal{F}_{n}$ containing $x.$

Let $\beta _{1},$ $\beta _{2},$ $\gamma _{1}$ and $\gamma _{2}$ be real
numbers such that
\begin{equation*}
\frac{1}{2}<\beta _{1}<\gamma _{1}<\beta _{2}<\gamma _{2}<\frac{1}{3}.
\end{equation*}
We say that an interval $I_{j}\in \mathcal{F}_{n}$ is of type $1$
(respectively of type $2$) when
\begin{equation*}
\beta _{1}<\frac{N_{0}(j)}{n}<\gamma _{1}\text{ (respectively }\beta _{2}<%
\frac{N_{0}(j)}{n}<\gamma _{2}\text{).}
\end{equation*}
Let $I\in \mathcal{F}_{n}$ be of type $1$ (respectively of type $2$), put $%
\tilde{I}$ the set of intervals of order $n+6$ contained in $I$ and of the
same type that $I,$ also put $\check{I}$ the set of intervals of order $2n$
contained in $I$ and of type $2$ (respectively of type $1$).

Let $n_{0}\in \mathbb{N}$ be a multiple of $6$ and $\left( n_{p}\right) $ be
the sequence of integers defined by :
\begin{equation*}
n_{0},\text{ }n_{3i+1}=2^{n_{3i}}n_{0},\text{ }n_{3i+2}=2n_{3i+1}\text{ and }%
n_{3i+3}=2n_{3i+2}.
\end{equation*}
Remark that $n_{p}=n_{0}+6k,$ $k\in \mathbb{N}.$ \newline
For all $k\in \mathbb{N,}$ we construct the family $\mathcal{G}_{k}$ of
disjoined diadic intervals of order $n_{0}+6k$ such that $\mathcal{G}_{0}$
contains two intervals $I_{n_{0}}^{1}$ and $I_{n_{0}}^{2}$ respectively of
type $1$ and $2$, any element of $\mathcal{G}_{k+1}$ is contained in an
element of $\mathcal{G}_{k}$ that we call his father, all the elements of $%
\mathcal{G}_{k}$ give birth to the same number of son in $\mathcal{G}_{k+1}$
and to pass from $\mathcal{G}_{k}$ to $\mathcal{G}_{k+1}$ we distinguish the
three following cases :

\textbf{1}$^{\text{st}}$\textbf{\ case :} If $n_{3i}\leq n_{0}+6k<n_{3i+1},$
then for each $I\in \mathcal{G}_{k}$ we select two intervals in $\widetilde{I%
}.$ So $\mathcal{G}_{k+1}$ is the union of all these selected intervals.

\textbf{2}$^{\text{nd}}$\textbf{\ cas :} If $n_{3i+1}\leq n_{0}+6k<n_{3i+2},$
then for each $I$ $\in \mathcal{G}_{k}$ of type $1$ we select an interval in
$\widetilde{I}$, and for each $I\in \mathcal{G}_{k}$ of type $2$ we select
an interval $I_{j}$ of order $n_{0}+6(k+1)$ such that $\beta _{1}<\dfrac{%
N_{0}(j)}{n}<\gamma _{2}$ and containing at least an interval of order $%
n_{3i+2}$ and of type $1$. So $\mathcal{G}_{k+1}$ is the union of all these
selected intervals.

\textbf{3}$^{\text{rd}}$\textbf{\ cas :} If $n_{3i+2}\leq n_{0}+6k<n_{3i+3},$
then for each $I\in \mathcal{G}_{k}$ of type $1$ we select an interval in $%
\widetilde{I}$, and for each $I\in \mathcal{G}_{k}$ having an ancestor of
order $n_{3i+1}$ and of type $1$, we select an interval $I_{j}$ of order $%
n_{0}+6(k+1)$ such that $\beta _{1}<\dfrac{N_{0}(j)}{n}<\gamma _{2}$ and
containing at least an interval of order $n_{3i+3}$ and of type $2$. So $%
\mathcal{G}_{k+1}$ is the union of all these selected intervals.

Note that any\thinspace $I_{j}\in \left( \underset{k\geq 0}{\cup }\mathcal{G}%
_{k}\right) \mathcal{\ }$verifies $\beta _{1}<\dfrac{N_{0}(j)}{n}<\gamma
_{2}.$

An elementary calculus of counting assures us that the construction of the
family $\left( \underset{k\geq 0}{\cup }\mathcal{G}_{k}\right) $ is possible
for any $n_{0}$ big enough, also it permits us to impose the following
separation condition : \newline
for all $k\geq 0$, if $I,$ $J\in \mathcal{G}_{k}$ are of order $n$, then the
distance between $I$ and $J$ is bigger than $\dfrac{1}{2^{n-1}}$. Besides
for all $k\geq 1$, if $I\in \mathcal{G}_{k}$ is of order $n$, then the
distances between $I$ and his father's endpoints are bigger than $\dfrac{1}{%
2^{n}}.$

We associate the following relation on $\left( \underset{k\geq 0}{\cup }%
\mathcal{G}_{k}\right) $ : \newline
the two elements of $\mathcal{G}_{0}$ are in relation and two elements of $%
\mathcal{G}_{k+1}$ are in relation \ if their fathers, elements of $\mathcal{%
G}_{k}$, are in relation.

Thereafter we call selected interval any element of $\left( \underset{k\geq 0%
}{\cup }\mathcal{G}_{k}\right) .$

Put $p_{0}$, $p_{1}>0$ such that $p_{0}+p_{1}=1$ and let $\mu $ be a
probability measure on $\mathbb{R}$ such that
\begin{equation*}
\mu \left( \mathbb{R\backslash }\left[ 0,1\right[ \right) =0
\end{equation*}
and for all $I_{j}\in \mathcal{F}_{n}$ and $l\in \left\{ 0,1\right\} $,
\begin{equation*}
\mu \left( I_{jl}\right) =\left\{
\begin{array}{l}
p_{l}\text{ }\mu \left( I_{j}\right) \text{, if }I_{j}\text{ contains a
selected interval,} \\
\dfrac{\mu \left( I_{j}\right) }{2}\text{, otherwise.}
\end{array}
\right.
\end{equation*}
It is clear that supp$\mu =\left[ 0,1\right] $.

We first show that the infimum $s$ of $B_{\mu }$ is strictly negative, what
comes back to establish the following proposition.

\begin{proposition}
\label{prop4}
\begin{equation*}
\underset{q\rightarrow +\infty }{\lim }B_{\mu }(q)=-\infty .
\end{equation*}
\end{proposition}

\begin{proof}
First, let's remark that for all $I_{j}\in \mathcal{F}_{n}$,
\begin{equation}
p_{0}^{n}\leq \mu \left( I_{j}\right) \leq p_{1}^{n}.  \label{eq7}
\end{equation}
Let $(B(x_{i},r_{i}))_{i\in I}$ be a centered $\varepsilon -$packing of supp$%
\mu .$\newline For all $i\in I,$ let's consider the largest interval
$I_{n}(x_{i})$ included in $B(x_{i},r_{i}).$ It results that
$B(x_{i},r_{i})$ is covered by at the more two contiguous intervals
of $\mathcal{F}_{n-1}.$ It follows that
\begin{equation}
\frac{1}{2^{n}}\leq 2r_{i}\leq \frac{1}{2^{n-2}}  \label{eq8}
\end{equation}
and according to (\ref{eq7}), we obtain
\begin{equation}
p_{0}^{n}\leq \mu (B(x_{i},r_{i}))\leq 2p_{1}^{n-1}.  \label{eq9}
\end{equation}
From (\ref{eq8}), we deduce that for all $t\in \mathbb{R},$ there exist $%
c_{1}$, $c_{2}\in \mathbb{R}$ such that for all $n\in \mathbb{N}$,
\begin{equation}
\frac{c_{1}}{2^{nt}}\leq (2r_{i})^{t}\leq \frac{c_{2}}{2^{nt}},
\label{eq10}
\end{equation}
and from (\ref{eq9}), it follows that for all $q>0$,
\begin{equation}
p_{0}^{nq}\leq \mu \left( B\left( x_{i},r_{i}\right) \right)
^{q}\leq 2^{q}p_{1}^{(n-1)q}.  \label{eq11}
\end{equation}
Then, considering (\ref{eq10}) and (\ref{eq11}), there exists
$c_{3}\in \mathbb{R}$ such that
\begin{equation}
\mu \left( B\left( x_{i},r_{i}\right) \right) ^{q}(2r_{i})^{t}\leq
c_{3}2^{q}p_{1}^{(n-1)q}2^{-nt}.  \label{eq12}
\end{equation}
Otherwise, for all $n\in \mathbb{N}\backslash \left\{ 0\right\} $,
any
interval of $\mathcal{F}_{n-1}$, meets to the more two balls of $%
(B(x_{i},r_{i}))_{i\in I}$ verifying the relation
$\dfrac{1}{2^{n}}\leq 2r_{i}\leq \dfrac{1}{2^{n-2}}$, so according
to (\ref{eq12}), there exists a constant $C$ that only depends on
$q$ and $t$ such that
\begin{equation}
\underset{\frac{1}{2^{n}}\leq 2r_{i}\leq \frac{1}{2^{n-2}}}{\sum
}\mu \left( B\left( x_{i},r_{i}\right) \right) ^{q}(2r_{i})^{t}\leq
C(2p_{1}^{q}2^{-t})^{n}.  \label{eq13}
\end{equation}
For $\varepsilon >0$ small enough, while writing,
\begin{equation*}
\underset{i\in I}{\sum }\mu \left( B\left( x_{i},r_{i}\right)
\right) ^{q}(2r_{i})^{t}=\underset{n\geq 1}{\sum
}\underset{\frac{1}{2^{n}}\leq 2r_{i}\leq \frac{1}{2^{n-2}}}{\sum
}\mu \left( B\left( x_{i},r_{i}\right) \right) ^{q}(2r_{i})^{t},
\end{equation*}
it comes from the inequality (\ref{eq13}) that $\underset{i\in
I}{\sum }\mu
\left( B\left( x_{i},r_{i}\right) \right) ^{q}(2r_{i})^{t}<\infty $ while $%
t>1+q\dfrac{Logp_{1}}{Log2}.$ We deduce that
\begin{equation*}
\Lambda _{\mu }(q)\leq 1+q\frac{Logp_{1}}{Log2},
\end{equation*}
then, according to the Proposition \ref{prop1}\
\begin{equation*}
B_{\mu }(q)\leq 1+q\frac{Logp_{1}}{Log2},
\end{equation*}
finally
\begin{equation*}
\underset{q\rightarrow +\infty }{\lim }B_{\mu }(q)=-\infty .
\end{equation*}
\end{proof}

\begin{proposition}
\label{prop5}Put $B_{\mu -}^{^{\prime }}(1)$ the left derivative number of $%
B_{\mu }$ at $1$. Then
\begin{equation*}
B_{\mu -}^{^{\prime }}(1)\leq -1.
\end{equation*}
\end{proposition}

\begin{proof}
Let's recall that $B_{\mu }(1)=0$ and $B_{\mu }$ is convex. So to
prove that $B_{\mu -}^{^{\prime }}(1)\leq -1$, it is sufficient to
establish that for
all $q<1,$ $B_{\mu }(q)\geq 1-q,$what comes back to show that, according to (%
\ref{eq2}), if $\left( \underset{i}{\cup }E_{i}\right) $ is a
partition of supp$\mu $, then $\underset{i\in I}{\sum
}\overline{P}_{\mu }^{q,t}(E_{i})=\infty $.

Let's consider the case where for all $i\in I,$ $\overline{P}_{\mu
}^{q,t}(E_{i})<\infty $, the contrary case is obvious. Put $0<\varepsilon <%
\frac{1}{2^{n_{0}}}.$ For all $i\in I$, choose $\delta
_{i}<\varepsilon $ such that
\begin{equation}
\overline{P}_{\mu ,\delta _{i}}^{q,t}(E_{i})\leq \overline{P}_{\mu
}^{q,t}(E_{i})+\frac{1}{2^{i}}.  \label{eq14}
\end{equation}
According to the Besicovitch covering theorem \cite{Gu}, there
exists an integer $\zeta $ (that only depends on $\mathbb{R}$) such
that each $E_{i}$
is covered by $\overset{\zeta }{\underset{u=1}{\cup }}\left( \underset{j}{%
\cup }B\left( x_{ij},\delta _{i}\right) \right) $ and for all $1\leq
u\leq \zeta $, $\left( B\left( x_{ij},\delta _{i}\right) \right)
_{j}$ is a packing. Considering (\ref{eq14}), it follows that
\begin{equation*}
\underset{u=1}{\overset{\zeta }{\sum }}\underset{j}{\sum }\mu \left(
B\left( x_{ij},\delta _{i}\right) \right) ^{q}\left( 2\delta
_{i}\right) ^{t}\leq \zeta \left( \overline{P}_{\mu
}^{q,t}(E_{i})+\frac{1}{2^{i}}\right) .
\end{equation*}
Then,
\begin{equation}
\underset{i}{\sum }\left( \underset{u=1}{\overset{\zeta }{\sum }}\underset{j%
}{\sum }\mu \left( B\left( x_{ij},\delta _{i}\right) \right)
^{q}\left(
2\delta _{i}\right) ^{t}\right) \leq \zeta \underset{i}{\sum }\overline{P}%
_{\mu }^{q,t}(E_{i})+\zeta .  \label{eq15}
\end{equation}
Let's consider the sum
\begin{equation}
\underset{i}{\sum }\left( \underset{u=1}{\overset{\zeta }{\sum }}\underset{j%
}{\sum }^{^{\prime }}\mu \left( B\left( x_{ij},\delta _{i}\right)
\right) ^{q}\left( 2\delta _{i}\right) ^{t}\right)   \label{eq16}
\end{equation}
where $\underset{j}{\sum }^{^{\prime }}$ is taken on all $j$ such
that the distance between $x_{ij}$ and $I_{n_{0}}^{1}$ (respectively
$I_{n_{0}}^{2}$) is bigger than $\dfrac{1}{2^{n_{0}}}.$ In this
case, there exists $C\in \mathbb{R}$ that only depends on $n_{0}$
such that
\begin{equation*}
\mu \left( B\left( x_{ij},\delta _{i}\right) \right) \leq C\,m\left(
B\left( x_{ij},\delta _{i}\right) \right)
\end{equation*}
where $m$ is the Lebesgue measure. We deduce that
\begin{equation}
C^{q-1}\left( 2\delta _{i}\right) ^{q-1+t}\leq \mu \left( B\left(
x_{ij},\delta _{i}\right) \right) ^{q-1}\left( 2\delta _{i}\right)
^{t}\mu \left( B\left( x_{ij},\delta _{i}\right) \right) .
\label{eq17}
\end{equation}
Otherwise, the union of the balls that appear in the sum
(\ref{eq16}) recovers supp$\mu $ deprived of $I_{n_{0}}^{1},$
$I_{n_{0}}^{2}$ and the intervals of order $n_{0}$ that their are
contiguous. Therefore, according to (\ref{eq17}), we obtain
\begin{equation*}
\left( 1-\frac{6}{2^{n_{0}}}\right) C^{q-1}\left( 2\varepsilon
\right)
^{q-1+t}\leq \underset{i}{\sum }\left( \underset{u=1}{\overset{\zeta }{\sum }%
}\underset{j}{\sum }^{^{\prime }}\mu \left( B\left( x_{ij},\delta
_{i}\right) \right) ^{q}\left( 2\delta _{i}\right) ^{t}\right) .
\end{equation*}
We deduce that, while considering (\ref{eq15}),
\begin{equation*}
\left( 1-\frac{6}{2^{n_{0}}}\right) C^{q-1}\left( 2\varepsilon
\right) ^{q-1+t}\leq \zeta \underset{i}{\sum }\overline{P}_{\mu
}^{q,t}(E_{i})+\zeta .
\end{equation*}
Letting $\varepsilon \rightarrow 0,$ it results that $\underset{i\in
I}{\sum }\overline{P}_{\mu }^{q,t}(E_{i})=\infty $ while $t<1-q.$
What permits to establish that $B_{\mu }(q)\geq 1-q.$
\end{proof}

Consider the Cantor set
\begin{equation*}
\mathcal{C}=\underset{k\geq 1}{\cap }\left( \underset{I_{j}\in \mathcal{G}%
_{k}}{\cup }I_{j}\right)
\end{equation*}
and the function $g$ defined on $\left[ 0,1\right] $ by
\begin{equation*}
g(x)=-\frac{xLog\left( \frac{p_{0}}{p_{1}}\right) +Logp_{1}}{Log2}.
\end{equation*}

\begin{proposition}
\label{prop6} \ \ \ \ \ \newline
i. If $x\notin \mathcal{C},$ then
\begin{equation*}
\underset{r\rightarrow 0}{\lim }\frac{Log(\mu \left( B(x,r)\right) }{Log2r}%
=1.
\end{equation*}
ii. If $x\in \mathcal{C},$ then
\begin{equation*}
g(\beta _{1})\leq \underset{r\rightarrow 0}{\lim \inf }\frac{Log(\mu \left(
B(x,r)\right) }{Log2r}\leq \underset{r\rightarrow 0}{\lim \sup }\frac{%
Log(\mu \left( B(x,r)\right) }{Log2r}\leq g(\gamma _{2}).
\end{equation*}
\end{proposition}

\begin{proof}
i. Put $x\notin \mathcal{C}.$ Thanks to the separation condition,
for $r>0$ small enough, the ball $B(x,r)$ is contained in the union
of two contiguous
intervals of order $N,$ $I_{N}^{1}$ and $I_{N}^{2}$ that don't meet $%
\mathcal{C}$. For all interval of order $n,$ $I_{n}\subset
I_{N}^{1}\cup I_{N}^{2}$ there exist $c$, $c^{\prime }\in
\mathbb{R}$ such that
\begin{equation*}
\frac{c}{2^{n}}\leq \mu \left( I_{n}(x)\right) \leq \frac{c^{\prime }}{2^{n}}%
.
\end{equation*}
We deduce that
\begin{equation}
\underset{n\rightarrow +\infty }{\lim }\frac{Log(\mu \left(
I_{n}(x)\right) }{Log\left( \dfrac{1}{2^{n}}\right) }=1.
\label{eq18}
\end{equation}

Consider the largest interval $I_{n}(x)$ contained in the ball
$B(x,r)$, it follows that $B(x,r)$ is contained in the union of two
contiguous intervals of order $n-1,$ $I_{n-1}(x)$ and $J_{n-1}$,
thus
\begin{equation*}
\frac{c}{2^{n}}\leq \mu \left( B(x,r)\right) \leq \frac{c^{\prime
}}{2^{n}}
\end{equation*}
and
\begin{equation*}
\left| I_{n}(x)\right| \leq 2r\leq 2\left| I_{n-1}(x)\right| .
\end{equation*}
Therefore, from (\ref{eq18}), we obtain
\begin{equation*}
\underset{r\rightarrow 0}{\lim }\frac{Log(\mu \left( B(x,r)\right) }{Log2r}%
=1.
\end{equation*}

ii. It is clear that if $I_{j}\in \mathcal{G}_{k}$ is of order $n$,
then
\begin{equation*}
\mu \left( I_{j}\right) =p_{0}^{N_{0}(j)}p_{1}^{n-N_{0}(j)},
\end{equation*}
thus
\begin{equation}
\mu (I_{j})=\left| I_{j}\right| ^{g\left( \frac{N_{0}(j)}{n}\right)
}. \label{eq19}
\end{equation}
Otherwise, let's recall that
\begin{equation*}
\beta _{1}<\dfrac{N_{0}(j)}{n}<\gamma _{2}.
\end{equation*}
Since the function $g$ is strictly increasing, it follows that
\begin{equation}
g\left( \beta _{1}\right) <\dfrac{Log(\mu \left( I_{j}(x)\right)
}{Log\left| I_{j}(x)\right| }<g\left( \gamma _{2}\right) .
\label{eq20}
\end{equation}

Put $x\in \mathcal{C}$ and $r<\dfrac{1}{2^{n_{0}+6}}$. Thanks to the
separation condition, $B(x,r)$ is contained in one of the intervals $%
I_{n_{0}}^{1}$ or $I_{n_{0}}^{2}.$\newline Consider the smallest
interval $I_{n}(x)$ containing the ball $B(x,r)$, it follows, from
the separation condition, that if $B(x,r)$ doesn't contain the
selected interval $I_{n+6}(x)$, then it necessarily contains the
selected interval $I_{n+12}(x)$, therefore, we can write
\begin{equation*}
\mu \left( I_{n+12}(x)\right) \leq \mu \left( B(x,r)\right) \leq \mu
\left( I_{n}(x)\right)
\end{equation*}
and
\begin{equation*}
\left| I_{n+12}(x)\right| \leq 2r\leq \left| I_{n}(x)\right| .
\end{equation*}
From (\ref{eq20}), it results that
\begin{equation*}
g(\beta _{1})\leq \underset{r\rightarrow 0}{\lim \inf }\frac{Log(\mu
\left(
B(x,r)\right) }{Log2r}\leq \underset{r\rightarrow 0}{\lim \sup }\frac{%
Log(\mu \left( B(x,r)\right) }{Log2r}\leq g(\gamma _{2}).
\end{equation*}
\end{proof}

We stand thereafter in the case where $g(\gamma _{2})<1,$ Even if we choose $%
p_{0}>\gamma _{2}.$ Thus, according to the Proposition \ref{prop6},
\begin{equation*}
\overline{X}^{g(\gamma _{2})}=\mathcal{C.}
\end{equation*}

In all what follows, we choose the real number $\alpha $ such that
\begin{equation*}
g(\gamma _{1})<\alpha \leq g(\gamma _{2})\text{ and\ }\overline{X}^{\alpha
}\neq \varnothing .
\end{equation*}

\begin{proposition}
\label{prop7} \
\begin{equation*}
T_{\mu }(\alpha )\leq g(\gamma _{1})<\alpha .
\end{equation*}
\end{proposition}

\begin{proof}
Put $M\subset X_{\alpha }(\eta ,p)$ and $(B(x_{i},r_{i}))$ a centered $%
\varepsilon -$packing of $M.$\newline It is clear that for all $i\in
I$, $x_{i}\in \mathcal{C}.$ Then consider the largest selected
interval $I_{n}(x_{i})$ of order $n$, containing $x_{i}$ and
contained in $B(x_{i},r_{i}).$ It follows that
\begin{equation*}
\frac{1}{2^{n}}\leq 2r_{i}.
\end{equation*}
Consider the partition $I_{1}\cup I_{2}$ of $I$ such that\newline
\begin{equation*}
I_{1}=\left\{ i\in I:I_{n}(x_{i})\text{ is of type }1\right\} \text{ and }%
I_{2}=I\backslash I_{1}.
\end{equation*}
Let's recall that, any interval $I_{n}(x_{i})$, $i\in I_{2}$, is in
relation
with an unique selected interval of order $n$ and of type $1$ centered in $%
x_{i}^{\prime }\in M$ that is denoted by $I_{n}(x_{i}^{\prime }).$
Thanks to
the separation condition, $\left( B\left( x_{i}^{\prime },\dfrac{1}{2^{n}}%
\right) \right) _{i\in I_{2}}$ is a centered $\varepsilon -$packing
of $M.$ Then we consider the family $\left( B(y_{i},\delta
_{i})\right) _{i\in I}$ indexed by $I$ and defined by
\begin{equation*}
B(y_{i},\delta _{i})=\left\{
\begin{array}{l}
B(x_{i},r_{i})\text{, }i\in I_{1} \\
B\left( x_{i}^{\prime },\frac{1}{2^{n}}\right) ,\text{ }i\in I_{2}.
\end{array}
\right.
\end{equation*}
We verify that
\begin{equation*}
\frac{Log\mu (B(y_{i},\delta _{i}))}{Log2\delta _{i}}\leq
\frac{Log\mu
\left( I_{n}(x_{i})\right) }{Log\left( \dfrac{1}{2^{n}}\right) },\text{ }%
i\in I_{1}
\end{equation*}
and
\begin{equation*}
\frac{Log\mu (B(y_{i},\delta _{i}))}{Log2\delta _{i}}\leq
\frac{Log\mu
\left( I_{n}(x_{i}^{\prime })\right) }{Log\left( \dfrac{1}{2^{n}}\right) },%
\text{ }i\in I_{2}.
\end{equation*}
From (\ref{eq19}) and as $g$ is increasing, we deduce that for all
$i\in I,$
\begin{equation*}
\frac{Log\mu (B(y_{i},\delta _{i}))}{Log2\delta _{i}}\leq g(\gamma
_{1}).
\end{equation*}
Thus
\begin{equation*}
L_{\varepsilon ,(B(x_{i},r_{i}))_{i\in I}}^{2}(M)\leq g(\gamma
_{1}).
\end{equation*}
Then $L_{\varepsilon }^{2}(M)\leq g(\gamma _{1}),$ letting
$\varepsilon
\rightarrow 0$, we deduce that $L^{2}(M)\leq g(\gamma _{1}).$ The sequence $%
\left( L^{k}(M)\right) _{k}$ is decreasing, it follows that
\begin{equation*}
L(M)\leq g(\gamma _{1}),
\end{equation*}
therefore,
\begin{equation*}
T_{\mu }(\alpha )\leq g(\gamma _{1}),
\end{equation*}
but $g(\gamma _{1})<\alpha ,$ then
\begin{equation*}
T_{\mu }(\alpha )<\alpha .
\end{equation*}
\end{proof}

\begin{corollary}
\label{cor1}
\begin{equation*}
Dim(\overline{X}^{\alpha })\leq \frac{1}{\alpha }T_{\mu }(\alpha )\text{ }%
\underset{q\geq 1}{\inf }\left( \alpha q+B_{\mu }(q)\right) <\underset{q\geq
0}{\inf }(\alpha q+B_{\mu }(q)).
\end{equation*}
\end{corollary}

\begin{proof}
From the Proposition \ref{prop4} and the Theorem \ref{th2}, we
deduce the first inequality.

Otherwise, as $\alpha <1$ and from the Proposition \ref{prop5}, it
follows that $B_{\mu -}^{^{\prime }}(1)\leq -\alpha $, then
\begin{equation*}
\underset{q\geq 1}{\inf }\left( \alpha q+B_{\mu }(q)\right)
=\underset{q\geq 0}{\inf }\left( \alpha q+B_{\mu }(q)\right) .
\end{equation*}
therefore,
\begin{equation*}
\frac{1}{\alpha }T_{\mu }(\alpha )\text{ }\underset{q\geq 1}{\inf
}\left(
\alpha q+B_{\mu }(q)\right) =\frac{1}{\alpha }T_{\mu }(\alpha )\text{ }%
\underset{q\geq 0}{\inf }\left( \alpha q+B_{\mu }(q)\right) .
\end{equation*}
Finally, according to the Proposition \ref{prop7}, we deduce that
\begin{equation*}
\frac{1}{\alpha }T_{\mu }(\alpha )\text{ }\underset{q\geq 1}{\inf
}\left( \alpha q+B_{\mu }(q)\right) <\underset{q\geq 0}{\inf
}(\alpha q+B_{\mu }(q)).
\end{equation*}
\end{proof}

\end{document}